\documentclass[12pt]{article}
\usepackage{mathrsfs}
\usepackage{CJK}
\usepackage{bbm}
\usepackage{stmaryrd}
\usepackage{shortvrb,psfrag}
\usepackage{enumerate}
\usepackage{amsmath}
\usepackage{cases}
\usepackage{amsfonts}
\usepackage{latexsym,amssymb,amsmath,mathrsfs,amsbsy, amsthm}
\usepackage[dvips]{graphicx}
\usepackage[usenames]{color}
\usepackage{xspace,colortbl}
\usepackage{pdfsync}
\usepackage[dvipdfm,
pdfstartview=FitH, CJKbookmarks=true, bookmarksnumbered=true,
bookmarksopen=true,
citecolor=blue
,colorlinks=true
,pdfborder=1
,pdfstartpage=2
,pdfstartview=Fit]{hyperref}
\allowdisplaybreaks

\newtheorem{thm}{Theorem}[section]
\newtheorem{lem}[thm]{Lemma}

\newtheorem{rem}{Remark}

\newtheorem{prop}[thm]{Proposition}

\newcommand{\beq}{\begin{equation}}
\newcommand{\eeq}{\end{equation}}
\newcommand{\ben}{\begin{eqnarray}}
\newcommand{\een}{\end{eqnarray}}
\newcommand{\beno}{\begin{eqnarray*}}
\newcommand{\eeno}{\end{eqnarray*}}

\begin{document}
\title{{\bf Remarks on Liouville type theorems for the 3D steady axially symmetric Navier-Stokes equations
}
\\Wendong Wang
\\[2mm]
{\small Dalian University of Technology, China}\\
{\small \&University of Oxford, UK}\\[2mm]}

\date{\today}
\maketitle

\begin{abstract}
In this note, we investigate the 3D steady axially symmetric Navier-Stokes equations, and obtained Liouville type theorems if the velocity or the vorticity satisfies some a priori decay assumptions.
\end{abstract}

{\small {\bf Keywords:} Liouville type theorem, Navier-Stokes equations, axially symmetric Navier-Stokes equations}

\setcounter{equation}{0}
\section{Introduction}

An interesting question about Liouville type theorem of
the 3D stationary Navier-Stokes equations in $R^3$ is as follows: whether the solution of
\begin{align} \label{eq:SNS}\,\, \left\{
\begin{aligned}
&-\Delta u+u\cdot \nabla u=-\nabla{p},\\
&\quad\nabla\cdot u=0,\\
\end{aligned}
\right. \end{align}
satisfying the vanishing property at infinity
\ben\label{eq:vanishing at infty}
\lim_{|x|\rightarrow\infty}u(x)=0,
\een
and the bounded Dirichlet energy
\ben\label{eq:energy bound}
D(u)=\int_{R^3}|\nabla u|^2dx<\infty
\een
implies $u\equiv0$ is still an open problem, which is related to J. Leray (see also P12. Galdi \cite{Galdi}).

Many conditional criteria have been obtained for this issue. For example,
Galdi proved the above Liouville type theorem by assuming $u\in L^{\frac92}(R^3)$ in  \cite{Galdi}.
Chae in \cite{Chae} showed the condition $\triangle u\in L^{\frac65}(R^3)$ is sufficient for the vanishing property of $u.$ Also, Chae-Wolf gave an improvement of logarithmic form for Galdi's result in \cite{ChaeWolf}
by assuming that $\int_{R^3} |u|^{\frac92}\{\ln(2+\frac{1}{|u|})\}^{-1}dx<\infty$.
 Seregin obtained the conditional criterion $u\in BMO^{-1}(R^3)$ in \cite{Se}. Moreover, Kozonoa-Terasawab-Wakasugib proved $u\equiv0$ if the vorticity $w=o(|x|^{-\frac53})$ or
$\|u\|_{L^{\frac92,\infty}(R^3)}\leq \delta D(u)^{1/3}$ for a small constant $\delta$ in \cite{KTW}.  It is shown that all the above norms $u\in L^{\frac92}(R^3)$, the log form of $u\in L^{\frac92}(R^3)$ or $u\in L^{\frac92,\infty}(R^3)$ can be replaced by the norms in the annular domain $B_R\setminus B_{R/2}$ in \cite{SeWang} by Seregin and the author, where the following energy description was stated:
\beno
\int_{B_{R/2}}|\nabla u|^2dx&\leq&CR^{-2}\left(\int_{B_R\setminus B_{R/2}}|{u}|^2dx\right)+C(q)R^{2-\frac9q}\|{u}\|_{L^{q,\infty}(B_R\setminus B_{R/2})}^3
\eeno
where $B_R=B_R(0)$ is a ball centered at $0$ and $q>3.$ Note that the conditions (\ref{eq:vanishing at infty}) and (\ref{eq:energy bound}) are not used in \cite{SeWang} as in \cite{ChaeWolf}.
More references, we refer to \cite{ChaeYoneda,Se2,Se3} and the references therein.


Moreover,
the problem is not solved, even for the case of axially symmetric Navier-Stokes equations, to the best of the author's knowledge. Motivated by the result
Seregin in \cite{Se3}, where he proved that the condition $|u|\lesssim\frac{1}{|x'|^{\mu}}$ with $x'=(x_1,x_2)$ and $\mu\approx 0.77$ implies $u\equiv0$, we are aimed to improve the decay assumption. At first, let us introduce
the axially symmetric Navier-Stokes equations.
Let $u(x)=u_r(t,r,z)e_r+u_\theta(t,r,z)e_\theta+u_z(t,r,z)e_z $, where
\beno
&& e_r=(\frac{x_1}{r},\frac{x_2}{r},0)=(\cos\theta,\sin\theta,0),\\
&& e_\theta=(-\frac{x_2}{r},\frac{x_1}{r},0)=(-\sin\theta,\cos\theta,0),\\
&& e_z=(0,0,1)
\eeno
and (\ref{eq:SNS}) becomes
\begin{align} \label{eq:axi-NS}\,\, \left\{
\begin{aligned}
&b\cdot \nabla u_r-\Delta_0 u_r+\frac{u_r}{r^2}-\frac{u_\theta^2}{r}+\partial_r p=0,\\
&b\cdot \nabla u_\theta-\Delta_0 u_\theta+\frac{u_\theta}{r^2}+\frac{u_ru_\theta}{r}=0,\\
&b\cdot \nabla u_z-\Delta_0 u_z+\partial_z p=0,\\
&\partial_r(ru_r)+\partial_z(ru_z)=0,\\
\end{aligned}
\right. \end{align}
where
\beno
b=u_re_r+u_ze_z,\quad \Delta_0=\partial_{rr}+\frac1r\partial_r+\partial_{zz}.
\eeno
The vorticity is represented as
\beno
w=w_re_r+w_\theta e_\theta+w_z e_z=(-\partial_z u_\theta)e_r+(\partial_zu_r-\partial_r u_z)e_\theta+\frac{\partial_r(ru_\theta)}{r}e_z.
\eeno

There are also many developments on the Liouville type theorems of axi-symmetric case. For example,
Liouville type theorem was proved by assuming no swirl (i.e. $u_\theta=0$), see Koch-Nadirashvili-Seregin-Sverak \cite{KNSS} or  Korobkov-Pileckas-Russo\cite{KRR}. The condition
$r u_{\theta}\in L^q$ with some $q\geq 1$ or
$b\in L^3$ is enough,
see Chae-Weng in \cite{ChaeWeng}.
Specially, for the axially symmetric case, the decay of the velocity or the vorticity can be obtained: Choe-Jin \cite{CJ}, Weng \cite{Weng} proved that
\beno
&&|u_r(r,z)|+|u_z(r,z)|+|u_\theta(r,z)|\lesssim \sqrt{\frac{\ln r}{r}},\\
&&|w_\theta(r,z)|\lesssim r^{-(\frac{19}{16})^-},\quad |w_r(r,z)|+|w_z(r,z)|\lesssim r^{-(\frac{17}{16})^-}
\eeno
Recently, Carrillo-Pan-Zhang in \cite{CPZ} gave an alternative method for the decay of $u$ and an improvement for the decay bound of the vorticity
\beno
|w_\theta(r,z)|\lesssim {r^{-\frac54}}{(\ln r)^{\frac34}},\quad |w_r(r,z)|+|w_z(r,z)|\lesssim {r^{-\frac98}}{(\ln r)^{\frac{11}{8}}}
 \eeno
by using Brezis-Gallouet inequality.

It's a natural question: whether there exist the sharp constants $\mu_1,\mu_2$ such that $|(u_r(r,z),u_z(r,z),u_\theta(r,z))|\lesssim \frac{1}{r^{\mu_1}}$
or $|(w_r(r,z),w_z(r,z),w_\theta(r,z))|\lesssim \frac{1}{r^{\mu_2}}$ implies that $u\equiv0$ for the axially symmetric case?

With the help of energy estimates in \cite{SeWang} we can improve the result in \cite{Se3} to $\mu>\frac23$, which is almost a equivalent form of $u\in L^{\frac92,\infty}$.
\begin{thm}\label{thm:Liouville of axi-NS}
Suppose that $u$ is axially symmetric smooth solution of the equation (\ref{eq:axi-NS}) and for some $\mu >\frac23$,
  \beno
~~|u|\leq \frac{C}{(1+r)^{\mu}}.
\eeno
Then $u\equiv 0.$
\end{thm}

Note that $\Gamma=r u_{\theta}$ satisfies the special structure
\beno
b\cdot \nabla\Gamma-\triangle_0\Gamma+\frac{2}{r}\partial_r\Gamma=0
\eeno
and Maximum principle can be applied, thus the condition $u_\theta=o(\frac1r)$ as $|x|\rightarrow\infty$ implies $u$ is trivial. However, it's still known that whether $u_\theta=o(\frac1r)$ can be replaced by $u_\theta=O(\frac1r)$. But we show that the condition $|b|=O(\frac1r)$ or $b\in BMO^{-1}(R^3)$  is sufficient, which improved the assumption
$b\in L^3(R^3)$
 in \cite{ChaeWeng}.

Here we say a function $f\in BMO^{-1}(R^3)$ if there
exists a vector-value function $d\in R^{3}$ and  $d_{j} \in BMO(R^3)$ such that $f=div~ d=d_{j,j}$. It's well-known that for the BMO space, we have
\beno
\Gamma(s)=\sup_{x_0\in R^3,R>0}\left(\frac{1}{|B_R(x_0)|}\int_{B_R(x_0)}|d-d_{x_0,R}|^sdx\right)^{\frac1s}<\infty.
\eeno
for any $s\in [1,\infty)$.

In details, we obtained the following result.
\begin{thm}\label{thm:Liouville of axi-NS2}
Suppose that $u$ is axially symmetric smooth solution of the equation (\ref{eq:axi-NS}) satisfying (\ref{eq:vanishing at infty}) and (\ref{eq:energy bound}). Then $u\equiv0$ if
one of the following conditions is satisfied
\beno
&&
(i)~~b=(u_r,u_z)\in BMO^{-1}(R^3);\\
&&
(ii)~~|b|\leq \frac{C}{r}
\eeno
\end{thm}

For the decay of the vorticity, we also state the following corresponding result.
\begin{thm}\label{thm:Liouville of axi-NS3}
Suppose that $u$ is axially symmetric smooth solution of the equation (\ref{eq:axi-NS}) satisfying (\ref{eq:vanishing at infty}) and (\ref{eq:energy bound}). Moreover,
\beno
~~|(w_r, w_\theta, w_z)|\leq \frac{C}{r^{\beta}}, \quad\beta>\frac53.
\eeno
Then $u\equiv0$.
\end{thm}

\begin{rem}
This conclusion generalized the result of \cite{KTW} to the axially symmetric case, where  the condition $|w|=o({|x|^{-\frac53}})$ was put.
\end{rem}

Throughout this article, $C$ denotes a constant, which may be different from line to line. 


\section{Proof of Theorem \ref{thm:Liouville of axi-NS}}

Recall a Caccioppoli inequality in \cite{SeWang}, which is stated as follows.

\begin{prop}\label{lem:estimate of LQ}
Let $(u,p)$ be the smooth solution of (\ref{eq:SNS}). Then for $0<\delta\leq1$ and $\frac{6(3-\delta)}{6-\delta}<q<3$, we have
\beno\label{eq:increasing estimate-Lorentz2}
\int_{B_{R/2}}|\nabla u|^2dx&\leq&\frac{C}{R^2}\left(\int_{B_R\setminus B_{R/2}}|{u}|^2dx\right)\nonumber\\
&&+C(\delta)\left(\|{u}\|^{3-\delta}_{L^{q,\infty}(B_R\setminus B_{R/2})}R^{2-\frac{9-3\delta}{q}-\frac{\delta}{2}}\right)^{\frac{2}{2-\delta}}
\eeno
\end{prop}

%


{\bf Proof of Theorem \ref{thm:Liouville of axi-NS}.} Let  $C_R$ denote the cylindrical region $\{x; |x'|\leq R, |z|\leq R\}$, then it's easy to check that
\beno
B_R\subset C_R\subset B_{\sqrt{2}R}.
\eeno
Hence, it follows from Proposition \ref{lem:estimate of LQ} that
\ben\label{eq:increasing estimate-Lorentz2}
\int_{C_{\frac{\sqrt{2}}{4}R}}|\nabla u|^2dx&\leq&\frac{C}{R^2}\left(\int_{C_R\setminus C_{\frac{\sqrt{2}}{4}R}}|{u}|^2dx\right)\nonumber\\
&&+C(\delta)\left(\|{u}\|^{3-\delta}_{L^{q,\infty}(C_R\setminus C_{\frac{\sqrt{2}}{4}R})}R^{2-\frac{9-3\delta}{q}-\frac{\delta}{2}}\right)^{\frac{2}{2-\delta}}\nonumber\\
&\leq&C\|u\|_{L^{q}(C_R)}^2R^{1-\frac6q} +C(\delta,q)\left(\|{u}\|^{3-\delta}_{L^{q}(C_R)}R^{2-\frac{9-3\delta}{q}-\frac{\delta}{2}}\right)^{\frac{2}{2-\delta}}
\een
for $q>2$, where we used the property of Lorentz space
\beno
\|{u}\|_{L^{q,\infty}(\Omega)}\leq C(q,\ell) \|{u}\|_{L^{q,\ell}(\Omega)}
\eeno
(for example, see Proposition 1.4.10 in \cite{Gr}).

For $\mu q>2,$ we have
\beno
\|u\|_{L^{q}(C_R)}
\leq C\left(R\int_{0}^R(1+r)^{1-\mu q}dr\right)^{\frac1q}\leq C(\mu,q)R^{\frac1q}
\eeno
Then the terms of the right hand side of (\ref{eq:increasing estimate-Lorentz2}) is controlled by
\ben\label{eq:decay of cylinder}
\int_{C_{\frac{\sqrt{2}}{4}R}}|\nabla u|^2dx\leq C(\mu,q)R^{1-\frac4q}+C(\delta,\mu,q)\left(R^{2-\frac{\delta}{2}-\frac{6-2\delta}{q}}\right)^{\frac{2}{2-\delta}}
\een
{\bf Claim that:} for fixed $\mu>\frac23$,  there  exist constants $\delta\in (0,1)$ and $q$ such that
\ben\label{eq:less than zero}
\max\{6\frac{3-\delta}{6-\delta},\frac{2}{\mu}\}<q<3,\quad {\rm and}~~2-\frac{\delta}{2}-\frac{6-2\delta}{q}<0
\een
hence letting $R\rightarrow\infty$, by (\ref{eq:decay of cylinder}) we have
\beno
\int_{R^3}|\nabla u|^2dx=0,
\eeno
which implies $u\equiv 0.$

{\bf Proof of (\ref{eq:less than zero}).}
First for fixed $\mu>\frac23$, we choose $\delta_0\in(0,1)$ such that
\beno
\frac{2}{\mu}<4\frac{3-\delta_0}{4-\delta_0}
\eeno

Since $0<\delta_0<1$, we have
\beno
1-\frac{\delta_0}{4}<1-\frac{\delta_0}{6},
\eeno
and
\beno
6\frac{3-\delta_0}{6-\delta_0}<4\frac{3-\delta_0}{4-\delta_0}
\eeno
so we take
\beno
q=\frac12\left(\max\{6\frac{3-\delta_0}{6-\delta_0},\frac{2}{\mu}\}
+4\frac{3-\delta_0}{4-\delta_0}\right)
\eeno
Then
 we have
\beno
\max\{6\frac{3-\delta_0}{6-\delta_0},\frac{2}{\mu}\}<q<4\frac{3-\delta_0}{4-\delta_0}<3,
\eeno
which implies (\ref{eq:less than zero}).

Hence the proof of Theorem \ref{thm:Liouville of axi-NS} is complete.

\section{Proof of Theorem \ref{thm:Liouville of axi-NS2}}

Let $\phi(x)=\phi(r,z)\in C_0^\infty(C_R)$ and $0\leq \phi\leq 1$ satisfying
\begin{align*} \phi(x)=\left\{
\begin{aligned}
&1,\quad x\in C_{R/2},\\
&0, \quad x\in C_R^c
\end{aligned}
\right. \end{align*}
and
\beno
|\nabla\phi|\leq \frac{C}{R},\quad |\nabla^2\phi|\leq \frac{C}{R^2}.
\eeno
Without loss of generality, by Theorem X.5.1 in \cite{Galdi} we can assume that
\beno
\lim_{|x|\rightarrow\infty}|p|+|u|=0.
\eeno
Note that $\triangle p=-\partial_i\partial_j(u_iu_j)$, then using Calder\'{o}n-Zygmund estimates and gradient estimates of harmonic function, we have
\beno
\int_{R^3}|p|^3+|u|^6dx<CD(u)^3,
\eeno
and
\beno
\|\nabla p\|_{L^{\frac32}(R^3)}<CD(u),
\eeno
since $\||\nabla u|u\|_{ L^{\frac32}(R^3)}\leq CD(u)$.

Multiplying $\phi u\cdot$ on both sides of (\ref{eq:SNS}), integration by parts yields that
\beno
&&\int_{C_R}\phi\left(|\nabla u_r|^2+|\nabla u_\theta|^2+|\nabla u_z|^2+\frac{u_r^2}{r^2}+\frac{u_\theta^2}{r^2}\right)dx\\
&\leq &\int_{C_R}\left(\frac12|u|^2+p\right)(u_r\partial_r+u_z\partial_z) \phi dx+C\|u\|^2_{L^6(C_R\setminus C_{R/2})}\\
&\doteq& I+C\|u\|^2_{L^6(C_R\setminus C_{R/2})}
\eeno

{\bf Case (i).} Due to $u_r,u_z\in BMO^{-1}(R^3)$, we write
\beno
u_r=\partial_jd_{1,j},\quad u_z=\partial_jd_{2,j},\quad j=1,2,3,
\eeno
where $d_{1,j},d_{2,j}\in BMO(R^3)$. Also, denote $\bar{f}$ as the mean value of $f$ on the domain $C_R.$
Then  we have
\beno
I
&=&\int_{C_R}\left(\frac12|u|^2+p\right)\left[\partial_j(d_{1,j}-\bar{d}_{1,j})\partial_r+\partial_j(d_{2,j}-\bar{d}_{2,j})\partial_z\right] \phi dx\\
&=&-\int_{C_R}\partial_j\left(\frac12|u|^2+p\right)\left[(d_{1,j}-\bar{d}_{1,j})\partial_r\phi+(d_{2,j}-\bar{d}_{2,j})\partial_z\phi \right] dx\\
&&-\int_{C_R}\left(\frac12|u|^2+p\right)\left[(d_{1,j}-\bar{d}_{1,j})\partial_j(\partial_r\phi)+(d_{2,j}-\bar{d}_{2,j})\partial_j(\partial_z\phi)\right]  dx\\
\eeno
Recall that $\phi(x)=\phi(r,z)$ and
\beno
&&\partial_j\partial_z \phi= \partial_z\partial_j \phi,\quad {\rm for}~~j=1,2,3,\\
&&\partial_j\partial_r \phi= \partial_z\partial_j \phi,\quad {\rm for}~~j=3,\\
&& \partial_1\partial_r \phi= \cos\theta\partial_r^2 \phi,\quad \partial_2\partial_r \phi= \sin\theta\partial_r^2 \phi,
\eeno
which and the property of BMO function yield that
\beno
I&\leq &CR^{-1}\||\nabla(|u|^2)|+|\nabla p|\|_{L^{\frac32}(C_R\setminus{C_{R/2}})}(\|d_{1,j}-\bar{d}_{1,j}\|_{L^{3}(C_R)}+\|d_{2,j}-\bar{d}_{2,j}\|_{L^{3}(C_R)})\\
 &&+C R^{-2}(\|u\|^2_{L^6(C_R\setminus{C_{R/2}})}+\|p\|_{L^3(C_R\setminus{C_{R/2}})})
 (\|d_{1,j}-\bar{d}_{1,j}\|_{L^{\frac32}(C_R)}+\|d_{2,j}-\bar{d}_{2,j}\|_{L^{\frac32}(C_R)})\\
&\leq &C\||\nabla(|u|^2)|+|\nabla p|\|_{L^{\frac32}(C_R\setminus{C_{R/2}})}+C (\|u\|^2_{L^6(C_R\setminus{C_{R/2}})}+\|p\|_{L^3(C_R\setminus{C_{R/2}})})\\
&&
\rightarrow  0\quad (\rm{as}~~R\rightarrow\infty)
\eeno

Hence, the proof of case (i) is complete.

{\bf Case (ii).} When $|(u_r,u_z)|\leq \frac{C}{r}$ for  $r>0$,
\beno
I&=&\int_{C_R}\left(\frac12|u|^2+p\right)(u_r\partial_r+u_z\partial_z) \phi dx\\
&\leq &C\int_{C_R}\left(\frac12|u|^2+|p|\right)\left(\partial_r \ln(r)|\partial_r\phi|+ \partial_r \ln(r)|\partial_z\phi|\right)dx.\\
\eeno
Let $g(r)=\ln(r)$ and $\bar{g}$ be the mean value of $g$ on $\{x'; |x'|\leq R\}$. Then we have
\beno
I&\leq &-C\int_{C_R}\partial_r(\frac12|u|^2+|p|) (g-\bar{g})\left(|\partial_r\phi|+|\partial_z\phi|\right)dx\\
&&-C\int_{C_R}\left(\frac12|u|^2+|p|\right) (g-\bar{g})\partial_r\left(|\partial_r\phi|+|\partial_z\phi|\right)dx\\
&&-C\int_{C_R}\left(\frac12|u|^2+|p|\right) (g-\bar{g})\frac{1}{r}\left(|\partial_r\phi|+|\partial_z\phi|\right)dx\\
&\doteq& I_1+I_2+I_3
\eeno
Note that  $g\in BMO(R^2)$ (see, for example, Chapter IV \cite{Stein}), and we have
\beno
R^{-1}\left(\int_{C_R}|g-\bar{g}|^3dx\right)^{\frac13}\leq C\left(R^{-2}\int_{|x'|\leq R}|g-\bar{g}|^3dx\right)^{\frac13}\leq C
\eeno
and
\beno
R^{-2}\left(\int_{C_R}|g-\bar{g}|^{\frac23}dx\right)^{\frac23}\leq C,\quad R^{-3}\left(\int_{C_R}|g-\bar{g}|^{12}dx\right)\leq C
\eeno
Hence as the arguments of (i), we have
\beno
I_1+I_2\leq C\||\nabla(|u|^2)|+|\nabla p|\|_{L^{\frac32}(C_R\setminus{C_{R/2}})}+C (\|u\|^2_{L^6(C_R\setminus{C_{R/2}})}+\|p\|_{L^3(C_R\setminus{C_{R/2}})})
\eeno
For the term of $I_3$, we get
\beno
I_3&\leq& CR^{-1}(\|u\|^2_{L^6(C_R\setminus{C_{R/2}})}+\|p\|_{L^3(C_R\setminus{C_{R/2}})})\|g-\bar{g}\|_{L^{12}(C_R)}\|\frac1r\|_{L^{\frac{12}{7}}(C_R)}\\
&\leq&  CR^{-\frac14}(\|u\|^2_{L^6(C_R\setminus{C_{R/2}})}+\|p\|_{L^3(C_R\setminus{C_{R/2}})})\|g-\bar{g}\|_{L^{12}(C_R)}\\
&\leq&  C(\|u\|^2_{L^6(C_R\setminus{C_{R/2}})}+\|p\|_{L^3(C_R\setminus{C_{R/2}})})
\eeno

Hence, we can conclude that
\beno
I\rightarrow  0\quad (\rm{as}~~R\rightarrow\infty)
\eeno
The proof of Theorem \ref{thm:Liouville of axi-NS2} is complete.

\section{Proof of Theorem \ref{thm:Liouville of axi-NS3}}

We are going to prove that
\begin{prop}\label{lem:ur decay from w} Assume that the conditions of Theorem \ref{thm:Liouville of axi-NS3} hold. (1) Let $w_\theta\leq C r^{-\beta}$ with $\beta>1$. Then
we get for $r>1$
\beno
~~|u_r(r,z)|+|u_z(r,z)|\leq C\left\{
\begin{aligned}
& (1+r)^{-\frac32+\frac{1}{2(\beta-1)}},\quad \beta>2,\\
&(1+r)^{1-\beta},\quad 1<\beta<2,\\
& (1+r)^{-1}\ln(r+1),\quad \beta=2.
\end{aligned}
\right.
\eeno
(2) Let $|w_r|+|w_z|\leq C r^{-\beta}$ with $\beta>1$. Then
we get for  $r>1$
\beno
~~|u_\theta(r,z)|\leq C\left\{
\begin{aligned}
& (1+r)^{-\frac32+\frac{1}{2(\beta-1)}},\quad \beta>2,\\
&(1+r)^{1-\beta},\quad 1<\beta<2,\\
& (1+r)^{-1}\ln(r+1),\quad \beta=2.
\end{aligned}
\right.
\eeno
\end{prop}

{\bf Proof of Theorem \ref{thm:Liouville of axi-NS3}.} It follows from Proposition \ref{lem:ur decay from w} and Theorem \ref{thm:Liouville of axi-NS} directly.

Next we are aimed to prove Proposition \ref{lem:ur decay from w}.
Firstly, we introduce a representation formula of $u_r$, $u_z$ and $u_\theta$ with the help of the vorticity. Since $b=u_re_r+u_ze_z$ and
\beno
\nabla\times b=w_\theta e_\theta,\quad \nabla\times (u_\theta e_\theta)=w_r e_r+w_z e_z
\eeno
by Biot-Savart law, we can get the integral representation of the velocity as follows(for example, see Lemma 2.2 for a local version by
Choe-Jin \cite{CJ}, also see Lemma 3.10 by Weng \cite{Weng}).
\begin{lem} Like the vorticity at the point $(r\cos\theta,r\sin\theta, z)$ denoted by $(w_r,w_\theta,w_z)(r,z)$, we write the vorticity at the point $(\rho\cos\phi,\rho\sin\phi, k)$ as $(w_\rho,w_\phi,w_k)(\rho,k)$. Then we have
\ben\label{eq:ur}
u_r(r,z)&=&\int_{-\infty}^{\infty}\int_0^\infty \Gamma_1(r,\rho,z-k)w_\phi(\rho,k)\rho d\rho dk,
\een
\ben\label{eq:uz}
u_z(r,z)&=&-\int_{-\infty}^{\infty}\int_0^\infty \Gamma_2(r,\rho,z-k)w_\phi(\rho,k)\rho d\rho dk
\een
\ben\label{eq:u theta}
u_\theta(r,z)&=&\int_{-\infty}^{\infty}\int_0^\infty \Gamma_3(r,\rho,z-k)w_k(\rho,k)\rho d\rho dk\nonumber\\
&&\quad -\int_{-\infty}^{\infty}\int_0^\infty \Gamma_1(r,\rho,z-k)w_\rho(\rho,k)\rho d\rho dk
\een
where
\beno
&&\Gamma_1(r,\rho,z-k)=\frac{1}{4\pi}\int_0^{2\pi}\frac{z-k}{[r^2+\rho^2-2r\rho\cos\phi +(z-k)^2]^{\frac32}}\cos\phi d\phi\\
&&\Gamma_2(r,\rho,z-k)=-\frac{1}{4\pi}\int_0^{2\pi}\frac{\rho-r\cos\phi}{[r^2+\rho^2-2r\rho\cos\phi +(z-k)^2]^{\frac32}}d\phi,\\
&&\Gamma_3(r,\rho,z-k)=-\frac{1}{4\pi}\int_0^{2\pi}\frac{\rho-r\cos\phi}{[r^2+\rho^2-2r\rho\cos\phi +(z-k)^2]^{\frac32}}\cos\phi d\phi.
\eeno
\end{lem}

Secondly, we give the bounds of estimate of $\Gamma_2$, $\Gamma_3$ and $\Gamma_1$, which will be used in the proof. This is similar to that in \cite{CJ}, where $\rho\approx r $ was assumed. Here we consider all $\rho>0$ and large $r>0.$ In details, we have the following estimates.
\begin{lem}[Estimate of $\Gamma_2$, $\Gamma_3$ and $\Gamma_1$]\label{lem: bound of Gamma2}
\ben\label{eq:estimate of Gamma 2}
|\Gamma_2(r,\rho,z-k)|+|\Gamma_3(r,\rho,z-k)|\leq \frac{C}{(\max\{\rho,r\})^{\alpha}[(r-\rho)^2+(z-k)^2]^{\frac{2-\alpha}2}},
\een
for $r>1$ and $0\leq\alpha\leq 1$;
\ben\label{eq:estimate of Gamma 5}
|\Gamma_1(r,\rho,z-k)|\leq \frac{C|z-k|}{(\max\{\rho,r\})^{\alpha}[(r-\rho)^2+(z-k)^2]^{\frac{3-\alpha}2}},
\een
where $r>1$, $0\leq\alpha\leq 1$ for $\frac{r}{4}\leq \rho\leq 4r,$ and $0\leq\alpha\leq 3$ for $\rho<\frac{r}{4}$ or $\rho\geq 4r.$
\end{lem}

Thirdly, we assume Lemma \ref{lem: bound of Gamma2} holds and complete the proof of Proposition \ref{lem:ur decay from w} and Lemma \ref{lem: bound of Gamma2} is proved later.

{\bf Proof of Proposition \ref{lem:ur decay from w}:} At first, we estimate the term of $u_r(r,z)$.
Let
\beno
I&=&u_r(r,z)=\int_{-\infty}^{\infty}\int_0^\infty \Gamma_1w_\phi\rho d\rho dk\\
&=&\int_{-\infty}^{\infty}\int_0^{r^\gamma/8} \Gamma_1w_\phi\rho d\rho dk+\int_{-\infty}^{\infty}\int_{r^\gamma/8}^{r/4} \Gamma_1w_\phi\rho d\rho dk+\int_{-\infty}^{\infty}\int_{r/4}^{r-r^\delta/2}\Gamma_1w_\phi\rho d\rho dk\\
&&+\int_{-\infty}^{\infty}\int_{r-r^\delta/2}^{r+r^\delta/2}\Gamma_1w_\phi\rho d\rho dk+ \int_{-\infty}^{\infty}\int_{r+r^\delta/2}^{4r}\Gamma_1w_\phi\rho d\rho dk +\int_{-\infty}^{\infty}\int_{4r}^\infty\Gamma_1w_\phi\rho d\rho dk\\
&=&I_1+\cdots+I_6,
\eeno
where $0\leq\gamma,\delta\leq 1$, to be decided.

For the term $I_1$, by (\ref{eq:estimate of Gamma 5}) and $\|w_\phi\|_{L^2(R^3)}^2\leq C D(u)<\infty$ we get
\beno
I_1&\leq&C \left(\int_{-\infty}^{\infty}\int_0^{r^\gamma/8}|\Gamma_1(r,\rho,z-k)|^2\rho d\rho dk\right)^{\frac12}\\
&\leq&C \left(\int_{-\infty}^{\infty}\int_0^{r^\gamma/8}\frac{|z-k|^2}{r^{2\alpha}[r^2+(z-k)^2]^{3-\alpha}}\rho d\rho dk\right)^{\frac12}\\
&\leq&C r^{-\frac32}\left(\int_{-\infty}^{\infty}\int_0^{r^\gamma/8}\frac{r^{-2}|z-k|^2}{[1+r^{-2}(z-k)^2]^{3-\alpha}} r^{-1}dk~\rho d\rho\right)^{\frac12}\leq C r^{-\frac32+\gamma}
\eeno
where $0\leq \alpha<\frac32.$

For the term $I_2$, using $r>1$, (\ref{eq:estimate of Gamma 5}) and $w_\theta\leq C r^{-\beta}$ 
\beno
I_2&\leq&C \int_{-\infty}^{\infty}\int_{r^\gamma/8}^{r/4} \Gamma_1\rho^{1-\beta} d\rho dk\\
&\leq&C \left(\int_{-\infty}^{\infty}\int_{r^\gamma/8}^{r/4}\frac{|z-k|}{r^{\alpha}[r^2+(z-k)^2]^{\frac{3-\alpha}{2}}}\rho^{1-\beta} d\rho dk\right)\\
&\leq& C\left\{
\begin{aligned}
& r^{-1+\gamma(2-\beta)}\quad (\beta>2)\\
& r^{-1}\ln r\quad (\beta=2)\\
& r^{1-\beta}\quad (1<\beta<2)
\end{aligned}
\right.
\eeno
where $0\leq \alpha<1.$

Moreover, for the term $I_3$, by (\ref{eq:estimate of Gamma 5}) and $w_\theta\leq C r^{-\beta}$ 
\beno
I_3&\leq&C \int_{-\infty}^{\infty}\int_{r/4}^{r-r^\delta/2} \Gamma_1\rho^{1-\beta} d\rho dk\\
&\leq&C \left(\int_{-\infty}^{\infty}\int_{r/4}^{r-r^\delta/2}\frac{|z-k|}{r^{\alpha}[(r-\rho)^2+(z-k)^2]^{\frac{3-\alpha}{2}}}\rho^{1-\beta} d\rho dk\right)\\
&\leq&C r^{-\alpha-\delta+\alpha\delta}\left(\int_{-\infty}^{\infty}\int_{r/4}^{r-r^\delta/2}\frac{r^{-\delta}|z-k|}{[\frac14+r^{-2\delta}(z-k)^2]^{\frac{3-\alpha}{2}}} r^{-\delta}dk\rho^{1-\beta} d\rho\right)\\
&\leq&C\left\{
\begin{aligned}
& r^{2-\beta-\alpha-\delta+\delta\alpha}\quad (\beta<2~or~\beta>2)\\
& r^{-\alpha-\delta+\alpha\delta}\ln r\quad (\beta=2)
\end{aligned}
\right.
\eeno
where $0\leq \alpha<1$.

Similarly, for $I_5$ we have
\beno
I_5&\leq&C\left\{
\begin{aligned}
& r^{2-\beta-\alpha-\delta+\delta\alpha}\quad (\beta<2~or~\beta>2)\\
& r^{-\alpha-\delta+\alpha\delta}\ln r\quad (\beta=2)
\end{aligned}
\right.
\eeno
where $0\leq \alpha<1$.

Furthermore, for $0\leq\alpha<1$ by (\ref{eq:estimate of Gamma 5}) and $w_\theta\leq C r^{-\beta}$  we have
\beno
I_4&\leq&C \int_{-\infty}^{\infty}\int_{r-r^\delta/2}^{r+r^\delta/2} \Gamma_1\rho^{1-\beta} d\rho dk\\
&\leq&C \left(\int_{-\infty}^{\infty}\int_{r-r^\delta/2}^{r+r^\delta/2}\frac{|z-k|}{r^{\alpha}[(r-\rho)^2+(z-k)^2]^{\frac{3-\alpha}{2}}}\rho^{1-\beta} d\rho dk\right)\\
&\leq&C \left(\int_{r-r^\delta/2}^{r+r^\delta/2}r^{-\alpha}(r-\rho)^{-1+\alpha}\rho^{1-\beta} d\rho \right)\\
&\leq&Cr^{1-\beta-\alpha} \left(\int_{r-r^\delta/2}^{r+r^\delta/2}(r-\rho)^{-1+\alpha} d\rho \right)\\
&\leq& C r^{1-\beta-\alpha+\delta\alpha}\quad (\beta>1)
\eeno

Finally, (\ref{eq:estimate of Gamma 5}) and $w_\theta\leq C r^{-\beta}$ yield that
\beno
I_6&\leq&C \int_{-\infty}^{\infty}\int_{4r}^\infty\Gamma_1\rho^{1-\beta} d\rho dk\\
&\leq&C \left(\int_{-\infty}^{\infty}\int_{4r}^\infty\frac{|z-k|}{\rho^{\alpha}[\rho^2+(z-k)^2]^{\frac{3-\alpha}{2}}}\rho^{1-\beta} d\rho dk\right)\\
&\leq& C r^{1-\beta}\quad (\beta>1)
\eeno

Hence, concluding the estimates of $I_1,\cdots, I_6$, we have the following arguments.

{\bf Case a. $\beta>2$.} At this time, we have
\beno
I\leq C \left[r^{-\frac32+\gamma}+r^{-1+\gamma(2-\beta)}+r^{2-\beta-\alpha-\delta+\delta\alpha}+r^{1-\beta-\alpha+\delta\alpha}+ r^{1-\beta}\right]
\eeno
where $0\leq \alpha<1$ and $0\leq\gamma,\delta\leq 1$.

First, we choose $\gamma=\frac{1}{2(\beta-1)}$ such that $-\frac32+\gamma=-1+\gamma(2-\beta)$. Furthermore, we take $\alpha \uparrow 1,\delta\uparrow 1$ such that
\beno
(1-\delta)(1-\alpha)\leq \beta-\frac52+\frac{1}{2(\beta-1)}
\eeno
which implies
\beno
-1+\gamma(2-\beta)\geq 2-\beta-\alpha-\delta+\delta\alpha
\eeno
Moreover, note that
\beno
2-\beta-\alpha-\delta+\delta\alpha\geq 1-\beta\geq  1-\beta-\alpha+\delta\alpha
\eeno
Then, we get for $r>1$
\beno
|u_r(r,z)|\leq C r^{-\frac32+\frac{1}{2(\beta-1)}}.
\eeno

{\bf Case b. $\beta<2$.} At this time, we have
\beno
I\leq  C \left[r^{-\frac32+\gamma}+r^{2-\beta-\alpha-\delta+\delta\alpha}+r^{1-\beta-\alpha+\delta\alpha}+ r^{1-\beta}\right]
\eeno
where $0\leq \alpha<1$ and $0\leq\gamma,\delta\leq 1$. We choose $\gamma=0$ and  $\delta=1$, then we get
\beno
|u_r(r,z)|\leq C r^{1-\beta}.
\eeno

{\bf Case c. $\beta=2$.} At this time, we have
\beno
I\leq  C \left[r^{-\frac32+\gamma}+r^{-1}\ln r+ r^{-\alpha-\delta+\alpha\delta}\ln r+r^{1-\beta-\alpha+\delta\alpha}+ r^{1-\beta}\right]
\eeno
where $0\leq \alpha<1$ and $0\leq\gamma,\delta\leq 1$. We choose $\gamma=0$ and $\delta=1$, then we get
\beno
|u_r(r,z)|\leq C r^{-1}\ln r.
\eeno
Hence we complete the estimate of $u_r(r,z).$

Note that the bound of $\Gamma_1$ used as above is similar to the estimates of $\Gamma_2$ and $\Gamma_3$. Hence similar arguments hold for $u_z$ and $u_\theta$. The proof of Proposition \ref{lem:ur decay from w} is complete.

{\bf Proof of Lemma \ref{lem: bound of Gamma2}.}
The remaining part is devoted to proving Lemma \ref{lem: bound of Gamma2}, which is similar to that of \cite{CJ}, where the case $\frac{r}{4}<\rho<4r$ is discussed. Here we consider all the value $\rho>0$ and sketch the proof. First, for $k>0$ and $\beta\geq 1$ we find
\ben\label{eq:bound of an integral}
I=\int_{0}^{\frac{\pi}{2}}\frac{d\phi}{(\sqrt{1+k\sin^2\phi})^{\beta}}\leq
\left\{
\begin{aligned}
&C(\delta)\min\{1,k^{-\frac{\delta}{2}}\},\quad \beta=1\\
&C(\beta)\min\{1,k^{-\frac{1}{2}}\},\quad \beta>1
\end{aligned}
\right.
\een
for any $0\leq \delta<1$. Obviously, $k\leq C$ holds, and next we assume that $k$ is large enough.
Then for $0<\ell<1$
\beno
I\leq \ell+\int_{\ell}^{\frac{\pi}{2}}\frac{d\phi}{({k\sin^2\phi})^{\beta/2}}
\eeno
Due to $\phi\leq 2\sin\phi$ for $\phi\in (0,\frac{\pi}{2})$, we have
\beno
I\leq \ell+2k^{-\beta/2}(\ln(\frac{\pi}{2})-\ln \ell),\quad \beta=1,
\eeno
and
\beno
I\leq \ell+2^\beta k^{-\beta/2}\frac{(\frac{\pi}{2})^{1-\beta}-\ell^{1-\beta}}{1-\beta},\quad \beta>1,
\eeno
which yield the required bound (\ref{eq:bound of an integral}) by choosing a suitable $\ell.$

Obviously, from the formulas of  $\Gamma_2, \Gamma_3$ and $\Gamma_1$, we have
\ben\label{eq:first estimate of Gamma2}
&&|\Gamma_i(r,\rho,z-k)|\leq \frac{\rho+r}{[(r-\rho)^2+(z-k)^2]^{\frac32}},\quad i=2,3;
\een
\ben\label{eq:first estimate of Gamma1}
&&|\Gamma_1(r,\rho,z-k)|\leq \frac{|z-k|}{[(r-\rho)^2+(z-k)^2]^{\frac32}}
\een
for all $\rho>0$ and $r>0.$

Next we go on estimating $\Gamma_2$, $\Gamma_3$, and $\Gamma_1$ carefully, respectively.

{\bf Step I.} Noting the periodic and even property and variable transform for $\phi$, we also have
\beno
\Gamma_2&=&-\int_0^{2\pi}\frac{1}{4\pi}\frac{\rho-r\cos\phi}{[r^2+\rho^2-2r\rho\cos\phi +(z-k)^2]^{\frac32}}d\phi\\
&=&-\int_0^{\frac{\pi}{2}}\frac{1}{\pi}\frac{\rho-r\cos2\phi}{[r^2+\rho^2-2r\rho\cos2\phi +(z-k)^2]^{\frac32}}d\phi
\eeno
and
\beno
\Gamma_{2}
&=&-\int_0^{\frac{\pi}{2}}\frac{1}{2\pi}\frac{\rho^2-2r\rho\cos2\phi+r^2+\rho^2-r^2}{\rho[(r-\rho)^2+4r\rho\sin^2\phi +(z-k)^2]^{\frac32}}d\phi\\
&\leq& C\frac{1}{\rho\sqrt{(r-\rho)^2+(z-k)^2}}\int_0^{\pi/2}\frac{d\phi}{\sqrt{1+K\sin^2\phi}}\\
&&-\frac{1}{2\pi}\frac{1}{\rho[(r-\rho)^2+(z-k)^2]^{\frac32}}\int_0^{\pi/2}\frac{\rho^2-r^2}{(\sqrt{1+K\sin^2\phi})^3}d\phi\\
&\doteq& I_1+I_2
\eeno
where
\beno
K=\frac{4r\rho}{(r-\rho)^2+(z-k)^2}
\eeno

When $K\leq 1$, that is ${4r\rho}\leq {(r-\rho)^2+(z-k)^2}$, we have ${(r-\rho)^2+(z-k)^2}\geq \frac12 r^2$ for $\rho\leq \frac{r}{2}$ and ${(r-\rho)^2+(z-k)^2}\geq 2 r^2$ for $\frac{r}{2}\leq \rho\leq 4r$. Moreover, for $\rho\geq 4r$ we have
\beno
{(r-\rho)^2+(z-k)^2}\geq (\frac34\rho)^2\geq (\frac35(\rho+r))^2\geq \frac{9}{25}(\rho+r)^2
\eeno
Hence for $K\leq 1$ we have
\ben\label{eq:Gamma2 K1}
\Gamma_{2}\leq C\frac{1}{\rho\sqrt{(r-\rho)^2+(z-k)^2}}
\een

When $K>1$, by (\ref{eq:bound of an integral}) we have
\ben\label{eq:Gamma2 K11}
\Gamma_2&\leq& C(\delta)\frac{1}{\rho\sqrt{(r-\rho)^2+(z-k)^2}}\nonumber\\
&&\cdot\left[ \left(\frac{(r-\rho)^2+(z-k)^2}{4r\rho}\right)^{\frac{\delta}{2}}+\frac{|\rho^2-r^2|}{(r-\rho)^2+(z-k)^2}\left(\frac{(r-\rho)^2+(z-k)^2}{4r\rho}\right)^{\frac12} \right]
\een
where $0\leq \delta<1$.

{\bf Case a. } For $r>1$ and $\rho\leq \frac{r}{4}$ or $\rho>4r$, by (\ref{eq:first estimate of Gamma2}) we know the estimate (\ref{eq:estimate of Gamma 2}) holds.

{\bf Case b. } For $r>1$ and $ \frac{r}{4}\leq \rho\leq 4r$ with $K\leq 1$, by (\ref{eq:first estimate of Gamma2}) and (\ref{eq:Gamma2 K1}) we know the estimate (\ref{eq:estimate of Gamma 2}) holds.

{\bf Case c. } For $r>1$ and $ \frac{r}{4}\leq \rho\leq 4r$ with $K>>1$, by (\ref{eq:first estimate of Gamma2}) and (\ref{eq:Gamma2 K11}) we know the estimate (\ref{eq:estimate of Gamma 2}) holds by noting that $(r-\rho)^2+(z-k)^2\leq 16r^2$ and
\beno
\frac{|\rho^2-r^2|}{(r-\rho)^2+(z-k)^2}\left(\frac{(r-\rho)^2+(z-k)^2}{4r\rho}\right)^{\frac12}\leq \frac{\rho+r}{\sqrt{4r\rho}}\leq 5.
\eeno

Hence the proof of $\Gamma_2$ is complete.

{\bf Step II.} The term $\Gamma_2$ is similar and we omitted the details.

{\bf Step III.} The term $\Gamma_1$ is estimated as follows.
\beno
\Gamma_1(r,\rho,z-k)&=&\frac{1}{2\pi}\int_0^{\pi}\frac{z-k}{[r^2+\rho^2-2r\rho\cos\phi +(z-k)^2]^{\frac32}}\cos\phi d\phi\\
&=&\frac{1}{\pi}\int_0^{\frac{\pi}{2}}\frac{z-k}{[(r-\rho)^2+4r\rho\sin^2\phi +(z-k)^2]^{\frac32}}\cos2\phi d\phi\\
&\leq&C\frac{|z-k|}{[(r-\rho)^2+(z-k)^2]^{\frac32}}\int_0^{\pi/2}\frac{1}{(\sqrt{1+K\sin^2\phi})^3}d\phi\\
&\doteq& I'
\eeno
where
\beno
K=\frac{4r\rho}{(r-\rho)^2+(z-k)^2}
\eeno

When $K\leq 1$, i.e. ${4r\rho}\leq {(r-\rho)^2+(z-k)^2}$, we have ${(r-\rho)^2+(z-k)^2}\geq \frac12 r^2$ for $\rho\leq \frac{r}{2}$ and ${(r-\rho)^2+(z-k)^2}\geq 2 r^2$ for $\frac{r}{2}\leq \rho\leq 4r$. Moreover, for $\rho\geq 4r$ we have
\beno
{(r-\rho)^2+(z-k)^2}\geq (\frac34\rho)^2
\eeno
Hence for $K\leq 1$ we have
\beno
{(r-\rho)^2+(z-k)^2}\geq \frac12(\max\{r,\rho\})^2
\eeno
Using (\ref{eq:first estimate of Gamma1}), for $K\leq 1$ we get
\ben\label{eq:Gamma1 K1}
|\Gamma_1(r,\rho,z-k)|\leq \frac{C|z-k|}{(\max\{\rho,r\})^{\alpha}[(r-\rho)^2+(z-k)^2]^{\frac{3-\alpha}2}},
\een
where $0\leq\alpha\leq 3$.

When $K>1$, i.e. ${4r\rho}\geq {(r-\rho)^2+(z-k)^2}$, which implies $\rho>\frac18 r$, by (\ref{eq:bound of an integral}) we have
\beno
|\Gamma_1(r,\rho,z-k)|&\leq& \frac{C|z-k|}{[(r-\rho)^2+(z-k)^2]^{\frac32}}\left(\frac{(r-\rho)^2+(z-k)^2}{4r\rho}\right)^{\frac12}\nonumber\\
&\leq& \frac{C|z-k|}{\sqrt{r\rho}[(r-\rho)^2+(z-k)^2]}
\eeno
Thus for $\frac18 r<\rho<4r$, we have
\ben\label{eq:Gamma1 K11}
|\Gamma_1(r,\rho,z-k)|\leq \frac{C|z-k|}{(\max\{\rho,r\})^{\alpha}[(r-\rho)^2+(z-k)^2]^{\frac{3-\alpha}2}}
\een
where $0\leq \alpha\leq 1$. For $\rho\geq 4r$, by (\ref{eq:first estimate of Gamma1}) we also derive that
\ben\label{eq:Gamma1 K12}
|\Gamma_1(r,\rho,z-k)|\leq \frac{C|z-k|}{(\max\{\rho,r\})^{\alpha}[(r-\rho)^2+(z-k)^2]^{\frac{3-\alpha}2}}
\een
where $0\leq \alpha\leq 3$.

Concluding the estimates (\ref{eq:Gamma1 K1}), (\ref{eq:Gamma1 K11}) and (\ref{eq:Gamma1 K12}), we complete the proof of the inequality
(\ref{eq:estimate of Gamma 5}).

\noindent {\bf Acknowledgments.}
W. Wang was supported by NSFC under grant 11671067,
 "the Fundamental Research Funds for the Central Universities" and China Scholarship Council.

\end{document}